\documentclass[12pt]{article}
\usepackage[reqno]{amsmath}
\usepackage{amssymb, theorem, enumerate}
\newcommand\al{\alpha}
\newcommand\be{\beta}

\newcommand\la{\lambda}

\newcommand\sg{\sigma}

\newcommand\cA{{\mathcal A}}

\newcommand\cD{{\mathcal D}}

\newcommand\cP{{\mathcal P}}

\newcommand\cx{{\mathbb C}}% complexes

\newcommand\re{{\mathbb R}}%reals

\newcommand\opk[1]{\mathop{\hbox{\rm #1}}\nolimits}

\newcommand\inv{^{-1}}
\newcommand\sbs{\subseteq}

\newcommand\seq[3]{#1_{#2},\ldots,#1_{#3}}
\newcommand\supp{\opk{supp}}

\newcommand\one{{\bf1}}
\newcommand\rk{\opk{rk}}

\expandafter\ifx\csname pdfoptionalwaysusepdfpagebox\endcsname\relax\else
\fi

\theoremstyle{change}
{\theorembodyfont{\slshape}
\newtheorem{theorem}{Theorem.}[section]
\newtheorem{lemma}[theorem]{Lemma.}
\newtheorem{corollary}[theorem]{Corollary.}}
\theorembodyfont{\rmfamily}

\newcommand\lref[1]{Lemma~\ref{lem:#1}}
\newcommand\tref[1]{Theorem~\ref{thm:#1}}
\newcommand\cref[1]{Corollary~\ref{cor:#1}}
\newcommand\sref[1]{Section~\ref{sec:#1}}

\def\proof{\noindent{{\sl Proof. }}}
\def\sqr#1#2{{\vbox{\hrule height.#2pt
    \hbox{\vrule width.#2pt height#1pt \kern#1pt
        \vrule width.#2pt}\hrule height.#2pt}}}
\def\eqed{\quad \vrule height7.5pt width4.17pt depth0pt} 
\def\qed{%
    \ifmmode\eqno\eqed
    \else\nobreak\ \hfill\eqed\medbreak\fi}

\newcommand\gbin[2]{\genfrac{[}{]}{0pt}{}{#1}{#2}}
\newcommand\qkne[3]{{#1}K_{{#2}:{#3}}}

\newcommand\qkvk{\qkne{q}{v}{k}}

\newcommand\nine{\{1,\ldots,9\}}
\newcommand{\pmat}[1]{\begin{pmatrix} #1 \end{pmatrix}}
\newcommand{\pnth}{\cP(3^3)}

\title{Independent Sets in Association Schemes}
\author{C. D. Godsil\footnotemark[1] and M. W. Newman\footnotemark[1]\footnote{Research supported by NSERC.} \\ Department of Combinatorics and Optimization \\ University
of Waterloo \\ 200 University Ave. W.\\ Waterloo, N2L 3G1, CANADA }
\begin{document}
\maketitle

\begin{abstract}
Let $X$ be $k$-regular graph on $v$ vertices and let $\tau$ denote the least eigenvalue of its adjacency matrix $A(X)$.  If $\al(X)$ denotes the maximum size of an independent set in $X$, we have the following well known bound:
\[
\al(X) \le\frac{v}{1-\frac{k}{\tau}}.
\]
It is less well known that if equality holds here and $S$ is a maximum independent set in $X$ with characteristic vector $x$, then the vector
\[
x-\frac{|S|}{v}\one
\]
is an eigenvector for $A(X)$ with eigenvalue $\tau$.  In this paper we show how this can be used to characterise the maximal independent sets in certain classes of graphs.  As a corollary we show that a graph defined on the partitions of $\{1,\ldots,9\}$ with three cells of size three is a core.
\end{abstract}

\section{Introduction}

Let $\pnth$ be the graph whose vertices are the partitions of $\nine$ with three cells of size three, where two such partitions are adjacent if each triple in one partition contains one point from each triple of the other.  (We will say that two such partitions are \textsl{skew}.)  Elementary arguments yield that $\pnth$ is a regular graph with valency 36 and with 280 vertices.

In \cite{kmbs}, Meagher and Stevens ask whether this graph is a core; this calls for some definitions.  If $X$ and $Y$ are graphs, a map $f$ from $V(X)$ to $V(Y)$ is a \textsl{homomorphism} if, whenever $u$ and $v$ are adjacent vertices in $X$, their images $f(x)$ and $f(y)$ are adjacent in $Y$.  An \textsl{endomorphism} of $X$ is a homomorphism from $X$ to itself.  Any automorphism of $X$ is an endomorphism and we say that $X$ is a \textsl{core} if all endomorphisms of $X$ are automorphisms.  For background on these topics, see \cite[Chapter 6]{cggr}.  

In this paper we prove that $\pnth$ is a core.  The central part of our argument is a determination of the independent sets of maximal size.  The method we use has other applications; we present two of these.

\section{Association Schemes}

An association scheme is a set of regular graphs whose edge sets partition the edges of the complete graph in a particularly nice way.  The precise definition is longer, and is stated in terms of the adjacency matrices of the graphs.  An \textsl{association scheme} $\cA$  with $d$ classes is a set $\seq A0d$ of $01$-matrices such that
\begin{enumerate}[(a)]
\item
$A_0=I$.
\item
$\sum_{i=0}^d A_i =J$.
\item
$A_i^T\in\cA$.
\item
$A_iA_j$ lies in the span of $\cA$ (over $\re$).
\item
$A_iA_j =A_jA_i$
\end{enumerate}
Here $J$ denotes the matrix with all entries equal to 1.  If the matrices in $\cA$ are all symmetric, we say that the scheme is \textsl{symmetric}.  All the schemes in this paper will be symmetric. 
In this case the matrices $\seq A1d$ are adjacency matrices of graphs $\seq X1d$.  Since the matrices in $\cA$ commute and since $J=\sum A_i$, we see that $A_iJ=JA_i$ for all $i$, whence the graphs $X_i$ must be regular.  We denote the number of vertices of these graphs by $v$, and the valency of $X_i$ by $v_i$.

Let $\cx[\cA]$ denote the span of $\cA$ over $\cx$.  This space is closed under matrix multiplication (by (d) above), and is called the \textsl{Bose-Mesner algebra} of the scheme.  
Since this algebra is closed under transpose and complex conjugation, the matrices in it can be simultaneously diagonalised.  More precisely, if the matrices $A_i$ are $v\times v$, the theory shows that that are subspaces $\seq U0d$ whose sum is $\cx^v$ and such that $U_j$ is contained in an eigenspace of $A_i$ for all $i$ and $j$.  It is traditional to denote the eigenvalue of $A_i$ on $U_j$ by $p_i(j)$, and to call the matrix $P$ such that
\[
P_{i,j} =p_j(i).
\]
Thus the $i$-th column of $P$ gives the eigenvalues of $A_i$.  The dimension of $U_j$ is denoted by $m_j$.  The \textsl{modified matrix of eigenvalues} is obtained by replacing the first column of $P$ with the multiplicities $m_i$.

Now let $E_j$ denote the orthogonal projection onto $U_j$.  Then the matrices $E_j$ form a second basis for $\cx[\cA]$.  (The problem is to show that they lie in $\cx[\cA]$.)  We have the following useful result due to Seidel.  For a proof, see \cite[Thm.~12.6.1]{bluebk}.

\begin{theorem}
\label{thm:seidel}
Suppose the matrices $\seq A0d$ form an association scheme on $v$ vertices with projections $\seq E0d$.  Then for any vector $x$,
\[
\sum_i \frac{x^TA_ix}{vv_i}A_i =\sum_j\frac{x^TE_jx}{m_j}E_j.
\] 
\end{theorem}

\section{Independent Sets}
\label{sec:3by3}
We start by giving a bound on the size of an independent set.  For a proof see, for example \cite[Lemma~9.6.2]{cggr}.  For graphs from association schemes (which is the only case we need), it follows from work of Delsarte \cite{DelTh}.  We call this result the \textsl{ratio bound} for $\al(X)$.

\begin{theorem}
\label{thm:ratbnd}
Let $X$ be a regular graph on $v$ vertices with valency $k$ and least eigenvalue $\tau$.  Then
\[
\al(X) \le\frac{v}{1-\frac{k}{\tau}};
\]
if equality holds and $x$ is the characteristic vector of an independent set of maximum size, then $x-\frac{\al(X)}v\one$ is an eigenvector for $A(X)$ with eigenvalue $\tau$.\qed
\end{theorem}

We apply this bound to $\pnth$.  Fortunately the eigenvalues of $\pnth$ have already been determined.

If $\pi$ and $\sg$ are two partitions of $\nine$ with three cells of size three, then their meet $\pi\wedge\sg$ has $3$, $9$, $7$, $6$ or $5$ cells.  Thus we can form four graphs $\seq X14$ with these partitions as their vertices, where partitions $\pi$ and $\sg$ are adjacent in $X_1$, $X_2$, $X_3$ or $X_4$ according as the number of cells of $\pi\wedge\sg$ is $9$, $7$, $6$ or $5$.
Mathon and Rosa \cite{part33} show that these four graphs are the graphs of an association scheme with four classes  with modified matrix of eigenvalues:
\[
\left(
\begin{array}{c|cccc}
 1& 36& 162& 54& 27\\
 27&-12&  -6&  6& 11\\
 48& 8&-6&-9&6\\
 120&2&-6&6&-3\\
 84&-4&12&-6&-3
\end{array}
\right)
\]
The second column of this matrix gives the eigenvalues of $\pnth$.  (For information on the eigenvalues of association schemes see, for example, \cite[Chapter~2]{bcn}.)  We see that $\tau=-12$, and that its multiplicity as an eigenvalue is 27.  We conclude that
\[
\al(\pnth) \le \frac{280}{1+\frac{36}{12}} =70.
\]
If $1\le i<j\le9$, there are exactly 70 partitions of $\nine$ with three cells of size three in which $i$ and $j$ lie in the same cell.  These 70 partitions form an independent set in $\pnth$, which we denote by $S_{i,j}$.  We aim to show that any independent set in $\pnth$ is of the form $S_{i,j}$, for some $i$ and $j$.

\section{Eigenvectors}

Let $M$ be the $280\times36$ $01$-matrix with rows indexed by the $3\times3$ partitions, columns indexed by the 2-element subsets of $\{1,\ldots,9\}$, and with the entry corresponding to the partition $\pi$ and the 2-element subset $ij$ equal to 1 if and only if $ij$ is contained in a cell of $\pi$.  Thus each row of $M$ is the characteristic vector of the set of 2-element subsets of $\{1,\ldots,9\}$ that are contained in a cell of $\pi$.  It follows that
\[
\frac19M\one =\one.
\]
If $B$ is the incidence matrix for $K_9$, the complete graph on nine vertices, then
\[
MB^T =2J.
\]
Since the rows of $B$ are linearly independent, this implies that $\rk(M)\le28$.  The columns of $M$ are the characteristic vectors of the independent sets $S_{i,j}$.  We will prove that if $x$ is the characteristic vector of an independent set of size $70$, then $x$ lies in the column space of $M$. 

\begin{lemma}
The columns of matrix $M-\frac14J$ are eigenvectors for $A$ with eigenvalue $-12$; they span a space with dimension 27.
\end{lemma}

\proof
The first claim is an immediate consequence of \tref{ratbnd}.  So we need only determine $\rk(M)$.  Since $\frac19M\one =\one$, the vector $\one$ lies in the column space of $M$.  Since 
\[
\one^T\Bigl(M-\frac14J\Bigr) =0,
\]
it will be enough to prove that $\rk(M)=28$.

Recall that $\rk(M)=\rk(M^TM)$.  Let $L$ denote the adjacency matrix of the line graph $L(K_9)$.  We note that
\[
|S_{1,2}|=70,\ |S_{1,2}\cap S_{1,3}| =10,\quad |S_{1,2}\cap S_{3,4}| =20
\]
and hence
\[
M^TM =70I+10L+20(J-I-L) =50I+20J-10L.
\]
The spectrum of $L$ is
\[
-2^{(27)},\ 5^{(8)},\ 14^{(1)}
\]
Accordingly $0$ is an eigenvalue of $M^TM$ with multiplicity eight (the eigenspace is spanned by the eigenvectors of $L$ with eigenvalue $5$), and therefore 
\[
\rk(M^TM)=28.
\]
Hence $\rk(M)=28$.\qed

\begin{corollary}
\label{cor:indM}
If $x$ is the characteristic vector of an independent set of size 70 in $\pnth$, then $x$ lies in the column space of $M$.
\end{corollary}

\proof
Set $z=x-\frac14\one$.  By the lemma, $z$ lies in the column space of $M$.  Since $\one$ also lies in the column space of $M$, it follows that $x$ does as well.\qed

\section{Structure of Independent Sets}

As noted earlier, the matrix $A$ lies in an association scheme with four classes.  Suppose that $A_1$, $A_2$, $A_3$ and $A_4$ are the matrices corresponding to the last four columns of the matrix of eigenvalues shown in \sref{3by3}.  Then $A_1=A$ and
\begin{equation}\label{mma}
MM^T =9I+2A_2+3A_3+5A_4.
\end{equation}
(The number of 2-element subsets of the cells of $\pi\wedge\sg$ is determined by the number of cells of $\pi\wedge\sg$.)  If the matrix idempotents of the Bose-Mesner algebra of the association scheme are $\seq E04$, then
\begin{equation}\label{mme}
MM^T =630E_0+70E_1.
\end{equation}
Here $E_0=\frac1{280}J$ and $E_1$ represents orthogonal projection onto the column space of $M-\frac14J$.  (One way to prove the above identity is to note that $MM^T$ and $M^TM$ have the same non-zero eigenvalues with the same multiplicities; alternatively this can be read off the matrix of eigenvalues of the scheme.)

\begin{lemma}
\label{lem:seid}
Suppose $x$ is the characteristic vector of an independent subset $S$ of $\pnth$ with size 70.  
Then
\[
\sum_{i=0}^4\frac{x^TA_ix}{vv_i}A_i
	=\frac14I+\frac{1}{18}A_2+\frac1{12}A_3+\frac{5}{36}A_4\\
	=\frac{70}4E_0+\frac{70}{36}E_1.
\]
\end{lemma}

\proof
We use Seidel's identity, \tref{seidel}.  Let $x$ be the characteristic vector of an independent set with size 70.  Then
\[
x^TE_0x =\frac{70\times 70}{280} =\frac{70}4.
\]
We also have
\[
E_1\one=0,\quad E_1\Bigl(x-\frac14\one\Bigr) =x-\frac14\one
\]
and therefore 
\[
x^TE_1x =\Bigl(x-\frac14\one\Bigr)^T\Bigl(x-\frac14\one\Bigr) 
	=70-\frac{70}4-\frac{70}4+\frac{280}{16}
	=\frac{210}4.
\]
Next
\[
E_2x =E_3x =E_4x =0
\]
and so for our particular values, the right side of Seidel's identity is
\[
\frac{70}{4}E_0+\frac{70}{36}E_1.
\]
Referring to \eqref{mme}, we see this equals $\frac1{36}MM^T$.  Using \eqref{mma}, we then find that
\[
\frac{70}{4}E_0+\frac{70}{36}E_1 
	=\frac1{36}(9I+2A_2+3A_3+5A_4)
\]
and the lemma follows from this.\qed

The valencies $v_1,\ldots v_4$ appear in the first row of the modified matrix of eigenvalues, the multiplicities $m_i$ are in the first column.  

\begin{corollary}
Let $S$ be an independent set of size 70 on $\pnth$.  If $\pi\in S$, then there are 36 partitions $\sg$ in $S$ such that $|\pi\wedge\sg|=7$, a further 18 such that $|\pi\wedge\sg|=6$ and $15$ such that $|\pi\wedge\sg|=5$.  
\end{corollary}

\section{The Characterisation}
\label{sec:maxinds}

We prove that an independent set of size 70 in $\pnth$ must be one of the sets $S_{i,j}$.  

Let $S$ be some independent set of size $70$ in $\pnth$, and $z$ its characteristic vector.  By
\tref{ratbnd} we know that $z-\frac14\one$ lies in the eigenspace of $X$ with eigenvalue $-12$
and, from \cref{indM}, it follows that $z$ lies in the column space of $M$.  From \lref{seid} we see that $S$ must contain two vertices $a$ and $b$ whose meet has exactly seven cells---two pairs and five singletons.

We have $z_a=z_b=1$ and $z_t=0$ for all vertices $t$ that are adjacent to $a$ or $b$.  Let $M_1$ be the submatrix of $M$ composed of the rows corresponding to the vertices adjacent to $a$ or $b$.  There is a vector $h$ such that $z=Mh$, and we then have $M_1h=0$.  Let $N$ be a matrix whose columns are a basis for the null space of $M_1$, so $h$ lies in the column space of $N$ and $z$ lies in the column space of $MN$.  Let $C$ be the matrix formed by the non-zero columns of the reduced column-echelon form of $MN$.  Then there is a vector $y$ such that $z=Cy$.  Since $C$ is in reduced column form and has full column rank, it has a set of rows that form an identity matrix of order $\rk(C)\times\rk(C)$.  Since $z$ is a $01$-vector, this implies that $y$ is a $01$-vector.

We now actually compute $C$ using \texttt{Maple}, and find that its rank is six.  Therefore $y$ is one of $2^6$ $01$-vectors and, using \texttt{Maple} again, we find that only three of these choices for $y$ have the property that $Cy$ is a $01$-vector:  one of these is the zero vector, and the other two are the characteristic vectors of the two $S_{ij}$ that contain $a$ and $b$.

The second \texttt{Maple} computation can be somewhat elided, as we now describe. It suffices to identify a submatrix $C_0$ corresponding to a subset of the rows of $C$ such that whenever
$y$ and $C_0y$ are $01$-vectors then $Cy$ is the characteristic vector of one of the two $S_{ij}$ that contain
$a$ and $b$, or the zero vector. We made no serious attempt to find an ``optimal'' $C_0$, but a casual glance
at $C$ reveals examples with $m=3$.

\section{A Core}

We show that, if the sets $S_{i,j}$ are the only maximum independent sets in $\pnth$, then $\pnth$ is a core.

Let $X$ denote $\pnth$ and suppose $\psi:X\to X$ is a homomorphism.  By \cite[Lemma~7.5.4]{cggr}, the preimage $\psi\inv(S_{i,j})$ is a maximum independent set in $X$, and hence is equal to $S_{k,\ell}$ for some $k$ and $\ell$.  Each maximal independent set is determined by a 2-element subset of $\nine$.  The independent sets corresponding to disjoint pairs have exactly 20 partitions in common, the sets corresponding to distinct overlapping pairs have 10 elements in common.  It follows that $\psi\inv$ determines a map from the vertex set of $L(K_9)$ to itself.  We show first that this mapping is an endomorphism of $L(K_9)$.

To this end, consider the sets
\[
\psi\inv(S_{1,2}\cap S_{1,i}),\quad i\in\{3,\ldots,9\}.
\]
If $\al$ is a partition in $\psi\inv(S_{1,2}\cap S_{1,i})$, then $\psi(\al)$ contains the cell $\{1,2,i\}$.  It follows that the six preimages above are pairwise disjoint.
Since
\[
\psi\inv(S_{1,2}\cap S_{1,i}) =\psi\inv(S_{1,2})\cap\psi\inv(S_{1,i})
\]
each of the preimages has size 10 or 20.  However the union of these preimages is $\psi\inv(S_{1,2})$, with size 70.  We conclude that if $i\in\{3,\ldots,9\}$, then
\[
|\psi\inv(S_{1,2}\cap S_{1,i})| =10
\]
and so we see that $\psi\inv$ determines an endomorphism of $L(K_9)$.

We will show below that $L(K_9)$ is a core. It follows from this that $\psi\inv$ in fact determines an
automorphism of $L(K_9)$. Therefore
\[
|\psi\inv(S_{i,j}\cap S_{k,\ell})|
\]
is equal to 70, 10 or 20 according as the pairs $ij$ and $k\ell$ are equal, distinct and overlapping or disjoint.

To complete the argument, consider the partition $\al$ with cells
\[
\{1,2,3\},\ \{4,5,6\},\ \{7,8,9\}.
\]
Then
\[
\{\al\} =S_{1,2}\cap S_{1,3} \cap S_{4,5}\cap S_{4,6}
\]
and
\[
\psi\inv(\al) =\psi\inv(S_{1,2})\cap \psi\inv(S_{1,3})\cap 
				\psi\inv(S_{4,5})\cap \psi\inv(S_{4,6}).
\]
Now $\{12,13\}$ and $\{45,46\}$ are vertex-disjoint edges in $L(K_9)$; since $\psi\inv$ induces an automorphism of $L(K_9)$, their images under $\psi\inv$ are vertex-disjoint edges too.  Therefore $\psi\inv(\al)\ne\emptyset$, and we have shown that $\psi$ is onto.  We conclude that $\pnth$ is a core.

It remains to show that $L(K_9)$ is a core.

\begin{lemma}
  $L(K_{2m+1})$ is a core.
\end{lemma}

\proof
Let $Y$ be the core of $L(K_{2m+1})$. Then $Y$ is vertex- and arc-transitive, and $\omega(Y)=2m$. Since $Y$ is
an induced subgraph of $X$, it is a line graph and therefore it is the line graph of an edge-transitive graph
$W$ with at least one vertex of valency $2m$. If $W$ is not bipartite, it is vertex-transitive and hence
regular; therefore $W=K_{2m+1}$. If $W$ is bipartite it is semi-regular. If $W$ has two vertices of degree
$2m$ they must be adjacent and therefore, since $W$ is edge-transitive, each edge joins two vertices of degree
$2m$ and $W=K_{2m+1}$. The only remaining case is when $W$ has a unique vertex of degree $2m$, whence
$W=K_{1,2m}$ and $Y=K_{2m}$. But $|V(Y)|$ divides $|L(K_{2m+1})|$. Since $2m$ does not divide $m(2m+1)$, we
conclude that $L(K_{2m+1})$ is its own core.
\qed

We note that the edge set of $L(K_{2m})$ can be partitioned into $2m-1$ perfect matchings. It follows that the
core of $L(K_{2m})$ is $K_{2m-1}$.

\section{$q$-Kneser Graphs}

Let $V$ be a set of size $v$.  The vertices of the \textsl{Kneser graph} $K_{v:k}$ are the $k$-subsets of $V$ and two $k$-subsets are adjacent if they are disjoint.  The Erd\H{o}s-Ko-Rado theorem asserts that an independent set in $K_{v:k}$ has size at most $\binom{v-1}{k-1}$, and that an independent set of this size consists of the $k$-subsets that contain a fixed element of $V$.  

We consider an analogous class of graphs.  Let $V$ be a vector space of dimension $v$ over the field of order $q$.  The vertices of the \textsl{$q$-Kneser graph} are the $k$-subspaces of $V$, and two $k$-subspaces are adjacent if their intersection is the zero subspace.  We denote this graph by $\qkvk$.  Our aim is to determine the independent sets of maximum size in $\qkvk$.

Let $q$ be a fixed positive integer.  We define
\[
[n] :=\frac{q^n-1}{q-1}.
\]
We then define $[n]!$ inductively by setting $[0]!=1$ and
\[
[n+1]! := [n+1][n]!.
\]
With this in hand we define the \textsl{$q$-binomial coefficient} $\gbin{v}{k}$ by
\[
\gbin{v}{k} :=\frac{[v]!}{[k]![v-k]!}.
\]
One reason the $q$-binomial coefficient is important is that $\gbin{v}{k}$ is the number $k$-dimensional subspaces of a $v$-dimensional vector space over a field with $q$ elements.

Suppose $V$ is a vector space of dimension $v$ over a field with $q$ elements.  Let $W_{1,k}(v)$ be the 01-matrix with rows indexed by the 1-dimensional subspaces, columns by the $k$-dimensional subspaces of $V$ and with $ij$-entry equal to 1 if and only if the $i$-th subspace of dimension $1$ is contained in the $j$-th subspace of dimension $k$.  Thus $W_{1,k}(v)$ is a $[v]$ by $\gbin{v}{k}$ matrix.  It is a standard result that the rows of $W_{1,k}(v)$ are linearly independent.  (See, e.g., \cite{kantor}.)   The rows of $W_{1,k}(v)$ are characteristic vectors of independent sets, namely the $k$-spaces that contain a given 1-dimensional subspace.

\begin{lemma}
If $v\ge2k$, the maximum size of an independent set in $\qkvk$ is
\[
\gbin{v-1}{k-1}.
\]
If $S$ is an independent set of this size, its characteristic vector lies in the column space of $W_{1,k}(v)^T$.
\end{lemma}

\proof
The least eigenvalue $\tau$ of $\qkvk$ is 
\[
-q^{k(k-1)}\gbin{v-k-1}{k-1},
\]
and the corresponding eigenspace has dimension $[v]-1$. This can be read off the matrix of eigenvalues of the
Johnson scheme, which is given in~\cite{Del76}. The valency of $\qkvk$ is
\[
q^{k^2}\gbin{v-k}{k}
\]
and hence the ratio bound applied to $\qkvk$ yields that
\[
\al(\qkvk) \le \gbin{v-1}{k-1}.
\]

It follows that the rows of $W_{1,k}(v)$ are characteristic vectors of independent sets from $\qkvk$ with maximal size.  Now $\rk(W_{1,k}(v))=[v]$ and $\one^T$ lies in the row space of $W_{1,k}(v)$.  Therefore the subspace of the column space of $W_{1,k}(v)^T$ orthogonal to $\one$ is the eigenspace of $\qkvk$ with eigenvalue $\tau$, and so the lemma follows.\qed

The following result is due to Hsieh~\cite{Hsi75} (his proof does not cover the case $v=2k+1, q=2$). Our proof
is new, and considerably simpler. It is related to the methods of Frankl and Wilson~\cite{FW86}, but they do not
give a proof of the characterization.
\begin{theorem}
If $v>2k$, then an independent set of maximum size in $\qkvk$ consists of the $k$-subspaces of $V$ that contain a given 1-dimensional subspace of $V$.
\end{theorem}

\proof
Let $M$ denote $W_{1,k}(v)^T$ and let $z$ be the characteristic vector of an independent set $S$ in $\qkvk$ with maximal size.   Then $z$ lies in the column space of $M$, and therefore there is a vector $h$ such that
\[
z=Mh
\]
If $\al$ is a $k$-subspace in $S$ and $\be$ is $k$-subspace such that $\be\cap\al=0$, then $\be\notin S$ and $z_{\be}=0$.  Hence $M_{\be}h=0$ for all $k$-subspaces $\be$ skew to $\al$.  

Let $B$ be a complement to the $k$-subspace $\al$, and let $M_{B}$ be the submatrix of $M$ formed by the rows $M_{\be}$, where $\be\sbs B$.   The columns of $M_{B}$ corresponding to the 1-dimensional subspaces not in $B$ are zero; the non-zero columns form the matrix $W_{1,k}(v-k)^T$, and therefore they are linearly independent.  Since $M_Bh=0$, it follows that the entries of $h$ indexed by the 1-dimensional subspaces of $B$ must all be zero.  Our choice of the complement $B$ was arbitrary, and consequently the entries of $h$ indexed by a 1-dimensional subspace not in $\al$ are zero.  Therefore the support of $h$ consists of 1-dimensional subspaces of $\al$ and so the $k$-subspaces in $S$ have a 1-dimensional subspace in common.\qed

It is not difficult to modify the above argument to obtain a proof of the Erd\H{o}s-Ko-Rado theorem.  The resulting proof is similar in its underlying ideas to the original proof by Wilson in \cite{rmw-ekr}.

\section{The Witt Graph on 77 Vertices}

Let $\cD$ be the well-known $3$-$(22,6,1)$ design, with 77 blocks.  The \textsl{Witt graph} $W$ has the 77 blocks of this design as its vertices; two blocks are adjacent if they are disjoint.  There are 21 blocks on each point of the design, and thus we find 22 independent sets of size 21.
The Witt graph is strongly regular, with eigenvalues and multiplicities
\[
(-6)^{(21)},\quad 2^{(55)},\quad 16^{(1)}.
\]
The ratio bound yields
\[
\al(W) \le\frac{77}{1+\frac{16}{6}} =21.
\]
Thus the independent sets given above have maximum size.  We aim to show that these are all the independent sets of maximum size in the graph.

Let $M$ be the $77\times22$ matrix with the characteristic vectors of our 22 independent sets as its columns.  (So $M^T$ is the usual incidence matrix of the design.)  Each pair of distinct points lies in exactly five blocks and consequently
\[
M^TM =16I_{22}+5J_{22},
\]
which is invertible.  Therefore $\rk(M)=22$.  Each column of
\[
M-\frac{3}{11}J
\]
is an eigenvector for $W$, and the columns of this matrix span the eigenspace with eigenvalue $-6$.  We conclude that if $x$ is the characteristic vector of an independent set of size 21, it lies in the column space of $M$.

Suppose $x$ contains the block $\al$, that is, $x_\al=1$ and assume $x=Mh$.  If $\be$ is a block disjoint from $\al$, then $x_\be=0$.  Let $M_1$ be the submatrix of $M$ formed by the rows corresponding to the blocks disjoint from $\al$.  Then $M_1h=0$.  Now $M_1$ is $16\times22$
\[
M_1 =\pmat{0&N},
\]
where $N$ is $16\times16$.  We show that $N$ is the incidence matrix of a $2$-$(16,6,2)$-design, whence it follows that $N$ is invertible.

Let $\al$ be a block of the Witt design and suppose $x,y$ are two points not in $\al$.  There are five blocks that contain $x$ and $y$ and, since any two distinct blocks have at most two points in common, these five blocks partition the 20 points other than $x$ and $y$.  It follows that three of these five blocks meet $\al$ in two points, and the remaining pair of blocks are disjoint from $\al$.  Hence the blocks disjoint from $\al$ form a 2-design on 16 points with block size six and with $\la=2$.  From this we deduce that there are exactly 16 blocks disjoint from $\al$ and so $N$ is the incidence matrix of a symmetric design.

Therefore $N$ is invertible and so it follows that $\supp(h)\sbs\al$.  Therefore each block in $\supp(x)$ contains $\supp(h)$, and consequently these 21 blocks have a point in common.

Since the Witt graph is the graph induced by the vertices not adjacent to a given vertex in the Higman-Sims
graph, our result leads to a determination of all independent sets of size $22$ in the Higman-Sims graph. It
is well known that the maximal independent sets of the Witt graph are as shown above: see for example
\cite{fihaem}, where the maximal independent sets in the Higman-Sims graph are determined by computer.
However our proof is new and gives further evidence that our method has a wide range of application.

\bibliographystyle{plain}
\bibliography{part33}

\end{document}